\documentclass[12pt,letterpaper]{article}

 \usepackage{amsmath,vmargin}
 \usepackage{theorem}
 \usepackage{amssymb}
 \usepackage{makeidx}

\theoremstyle{plain}
\newtheorem{thm}{Theorem}

\theoremstyle{plain}

\theoremstyle{plain}
\newtheorem{lem}{Lemma}

\theoremstyle{plain}

\theoremstyle{plain}
\newtheorem{cor}{Corollary}

\newcommand\proof{{\noindent{  Proof : }}}

\newcommand\F{\mathbb{F}}
\newcommand\R{\mathbb{R}}
\newcommand\C{\mathbb{C}}

\newcommand\dans{\longrightarrow}
\newcommand\associe{\longmapsto}

\newcommand\un{\underline}

\newcommand\FQH{\mathcal{F}_{q}\left(H\right)}
\newcommand\FQ{\mathcal{F}_{q}\left(H_{\C}\right)}
\newcommand\FO{\mathcal{F}_{0}\left(H_{\C}\right)}
\newcommand\FT{\mathcal{F}_{T}\left(H_{\C}\right)}

\newcommand\GO{\Gamma_{0}(H_{\R})}
\newcommand\GT{\Gamma_{T}(H_{\R})}
\newcommand\GQt{\Gamma_{q}(H_{\R},U_{t})}
\newcommand\GQtr{\Gamma_{q,r}(H_{\R}^{\prime},U_{t})}
\newcommand\HC{H_{\C}}

\newcommand\Eproof{\hfill{{\large $\square$}}}

\def\Dual#1#2{\langle #1 , #2 \rangle}

\def\fonction#1#2#3#4#5{
\begin{equation*}
\begin{split}
#1 : #2 &\dans #3 \\
     #4 &\associe #5
\end{split}
\end{equation*}}

\setmarginsrb{2cm}{2.2cm}{2cm}{2.2cm}{1mm}{6mm}{0mm}{10mm}

\begin{document}
\pagestyle{empty}
\renewcommand{\thefootnote}{}
\begin{center}
{ \Large \bf Non injectivity of the q-deformed von Neumann algebra }\\
\vspace{4 cm}
Alexandre Nou\\
Universit\'e de Franche-Comt\'e - Besancon\\
U.F.R des Sciences et Techniques\\
D\'epartement de Math\'ematiques\\
16 route de Gray - 25030 Besancon Cedex\\
nou@math.univ-fcomte.fr
\end{center}
\vspace{4 cm}
{\bf Abstract}. { In this paper we prove that the von Neumann algebra
  generated by q-gaussians is not injective as soon as the dimension
  of the underlying Hilbert space is greater than 1. Our approach is
  based on a suitable vector valued Khintchine type inequality for Wick products. The same proof also works for the more
  general setting of a Yang-Baxter deformation. Our techniques can also be extended to the so called q-Araki-Woods
  von Neumann algebras recently introduced by Hiai. In this latter case, we obtain the non injectivity under some asssumption on the spectral set of the positive operator asociated with the deformation.}
\footnotetext[1]{AMS classification : 46L65, 46L54 \hfill
Keywords : injectivity, deformation, Yang-Baxter operator}
\newpage

\section{Introduction}
\pagestyle{plain}
Let $H_{\R}$ be a real Hilbert space and $\HC$ its complexification. Let
  $T$ be a Yang-Baxter operator on
  $\HC\otimes\HC$ with $\|T\|<1$. Let $\FT$ be the associated deformed
  Fock space and $\GT$ the von Neumann algebra generated by the
  corresponding deformed
  gaussian random variables, introduced by Bozejko and Speicher
  \cite{BS} (also see \cite{BKS}). In addition, we will assume that
  $T$ is tracial, i.e that the vacuum expectation is a trace on $\GT$ (cf
  \cite{BS}). Under these assumptions, it was proved in \cite{BS} that
  $\GT$ is not injective as soon as $\dim H_{\R}> \frac{16}{(1-q)^{2}}$,
  where $\|T\|= q$. Since then the problem whether $\GT$ is not
  injective as soon as $\dim H_{\R} \ge 2$ had been left open. We emphasize
  that this problem remained open even in the particular case of the
  q-deformation, that is when $T=q\sigma$, where $\sigma$ is the
  reflexion : $\sigma(\xi\otimes \eta)=\eta \otimes \xi$. Recall that
  the free von Neumann algebra $\GO$ (corresponding to $T=0$) is not
  injective as soon as $n=\dim H_{\R} \ge 2$, for $\GO$ is isomorphic to
  the free group von Neumann algebra $VN(\F_{n})$ (cf. \cite{VDN}). The
  main result of this paper solves the above problem.

To explain the idea of our proof we first recall the main ingredient
  of the proof of the non injectivity theorem in \cite{BS}. It is the
  following vector-valued non-commutative Khintchine inequality. Let
  $(e_{i})_{i \in I}$ be an orthonormal basis of $H_{\R}$. Let $K$ be a
  complex Hilbert space and $B(K)$ the space of all bounded operators
  on $K$. Then for any finitely supported family $(a_{i})_{i \in
  I}\subset B(K)$

\begin{eqnarray*} \max\left\{\|\sum_{i\in I}a_{i}^{*}a_{i}\|_{B(K)}^{\frac{1}{2}}, \|\sum_{i\in I}a_{i}a_{i}^{*}\|_{B(K)}^{\frac{1}{2}}
\right\}
 & \le  & \left\|\sum_{i\in I}a_{i}\otimes G(e_{i})\right\|\\
  & \le &
\frac{2}{\sqrt{1-q}}
 \max\left\{\|\sum_{i\in I}a_{i}^{*}a_{i}\|_{B(K)}^{\frac{1}{2}}, \|\sum_{i\in I}a_{i}a_{i}^{*}\|_{B(K)}^{\frac{1}{2}}
\right\}
\end{eqnarray*}
where $G(e)=a^{*}(e)+a(e)$ is the deformed gaussian variable
associated with a vector $e \in H_{\R}$. Using this Khintchine inequality
and the equivalence between the injectivity and the semi-discreteness,
one easily deduces the non-injectivity of $\GT$ as soon as $\dim H_{\R}>
\frac{16}{(1-q)^{2}}$.

The proof of our non-injectivity theorem follows the same pattern. We
will first need to extend the preceding vector-valued non-commutative
Khintchine inequality to Wick products. It is well known that for any
$\xi$, a finite linear combination of elementary tensors, there is a
unique operator $W(\xi) \in \GT$ such that $W(\xi)\Omega=\xi$. Instead
of the previous inequality, the main ingredient of our proof is the following. Let $n\ge
1$. Let $(\xi_{\un{i}})_{|\un{i}|=n}$ be an orthonormal basis of
$\HC^{\otimes n}$ and $(\alpha_{\un{i}}) \subset B(K)$ a finitely
supported family. Then

\begin{equation}
\max_{0\le k \le n}\lbrace\|\sum_{|\un{i}|=n} \alpha_{\un{i}}\otimes R_{n,k}^{*}\xi_{\un{i}}\|\rbrace
\le
\|\sum_{|\un{i}|=n}\alpha_{\un{i}}\otimes W(\xi_{\un{i}})\|
\le
(n+1)C_{q}\max_{0\le k \le n}\lbrace\|\sum_{|\un{i}|=n} \alpha_{\un{i}}\otimes R_{n,k}^{*}\xi_{\un{i}}\|\rbrace
\end{equation}
where the norms in the left and right handside have to be taken in
$B(K)\otimes_{\min}H_{c}^{\otimes n-k}\otimes_{h}H_{c}^{\otimes k}$
(see Theorem 1 below for the
precise statement). Inequality (1) is the vector-valued version of
Bozejko's ultracontractivity inequality proved in \cite{B2} and thus it
solves a problem posed in \cite{B2}. Using (1) and a careful analysis
on the norms of Wick products on a same level, we deduce our
non-injectivity result.

The plan of this paper is as follows. The first section is devoted to
necessary definitions and preliminaries on the deformation by a Yang-Baxter operator and the associated von Neumann algebra. In this
section, we also include a brief discussion on the simplest case, the
free case, i.e. when $T=0$. All our results and arguments become very
simple in this case, for instance, inequality (1) above is then easy
to state and prove. The proof of the non-injectivity of $\GO$ can be
done in just a few lines. The reason why we have decided to include
such a discussion on the free case is the fact that it already
contains the main idea for the general case. In the second section we will
establish (1) and prove the non-injectivity of $\GT$. The last section aims at proving the non-injectivity of the
Araki-Woods factors $\Gamma_{q}(H,U_{t})$ introduced by Hiai in
\cite{H}. Note that Hiai proved a non-injectivity result with a condition on
the dimension of the spectral sets of the positive generator of $U_{t}$, which is similar to
that of \cite{BS}. The problem is left open whether the dimension can
go down to 2. Although we cannot completely solve this,
our method permits to improve in some sense the criterion for
non-injectivity given in \cite{H}.

\section{Preliminaries}

Recall that the free Fock space associated with $H_{\R}$ is given by
$$\FO = \bigoplus_{n \ge 0} \HC^{\otimes n}$$
where $\HC^{\otimes 0}$ is by definition $\C\Omega$ with $\Omega$ a
unit vector called the vacuum.

A Yang-Baxter operator on $\HC\otimes\HC$ is a self-adjoint
contraction satisfying the following braid relation :
$$(I\otimes T)(T\otimes I)(I\otimes T)=(T\otimes I)(I\otimes
T)(T\otimes I)$$
For $n\ge 2$ and $1\le k\le n-1$ we define $T_{k}$ on $\HC^{\otimes n}$
by
$$T_{k}=I_{\HC^{k-1}}\otimes T\otimes I_{\HC^{n-k-1}}$$
Let $S_{n}$ be the group of permutations on a set of $n$ elements. A
function $\varphi$ is defined on $S_{n}$ by quasi-multiplicative
extension of :
$$\varphi(\pi_{k})=T_{k}$$
where $\pi_{k}=(k,k+1)$ is the transposition exchanging $k$ and $k+1$, $1\le k \le n-1$. The symmetrizator
$P_{T}^{(n)}$ is the following operator defined on $\HC^{\otimes n}$ by :
$$P_{T}^{(n)}=\sum_{\sigma \in S_{n}}\varphi(\sigma)$$
$P_{T}^{(n)}$ is a positive operator on $\HC^{\otimes n}$ for any
Yang-Baxter operator $T$ and is strictly positive if $T$ is strictly
contractive (cf. \cite{BS}). In the latter case we are allowed to
define a new scalar product on $\HC^{\otimes n}$  (for $n\ge 2$) by :
$$\Dual{\xi}{\eta}_{T}=\Dual{\xi}{P_{T}^{(n)}\eta}$$
The associated norm is denoted by $\|.\|_{T}$.
The deformed Fock space associated with $T$ is then defined by
$$\FT=\bigoplus_{n\ge 0}\HC^{\otimes n}$$
where $\HC^{\otimes n}$ is now equipped with our deformed scalar
product for $n\ge 2$. From now on we will only consider a strictly
contractive Yang-Baxter $T$ and $\|T\|\le q < 1$.

For $f\in H_{\R}$, $a^{*}(f)$ will denote the creation operator associated
with $f$, and $a(f)$ its adjoint with respect to the T-scalar
product :
$$a^{*}(f)(f_{1}\otimes\dots \otimes f_{n})=f\otimes f_{1}\otimes \dots
\otimes f_{n}$$
For $f \in H_{\R}$ the deformed gaussian is the following hermitian
operator :
$$G(f)=a^{*}(f)+a(f)$$

\medskip

Throughout this paper we are interested in $\GT$ which is the von
Neumann algebra generated by all gaussians $G(f)$ for $f \in H_{\R}$ :
$$\GT=\{G(f) : f\in H_{\R}\}^{\prime\prime}\subset B(\FT)$$
Let $(e_{i})_{i\in I}$ be an orthonormal basis of $H_{\R}$ and set
$$t_{ij}^{sr}=\Dual{e_{s}\otimes e_{r}}{T(e_{i}\otimes e_{j})}$$
Then the following deformed commutation relations hold :
$$a(e_{i})a^{*}(e_{j})-\sum_{r,\; s \in
  I}t_{js}^{ir}a^{*}(e_{r})(e_{s})=\delta_{i j}$$
Moreover if the following condition holds
$$\Dual{e_{s}\otimes e_{r}}{Te_{i}\otimes e_{j}}=\Dual{e_{r}\otimes
  e_{j}}{Te_{s}\otimes e_{i}}$$
which is equivalent to the cyclic condition :
$$t_{ij}^{sr}=t_{si}^{rj}$$
then the vacuum is cyclic and separating for $\GT$ and the vacuum
expectation is a faithful trace on $\GT$ that will be denoted by $\tau$. If this cyclic condition
holds we say that $T$ is tracial, and from now on we will always assume
that $T$ has this property.

We will denote by $\Gamma_{T}^{\infty}(H_{\R})$ the subspace
$\GT\Omega$ of $\FT$. Since $\Omega$ is separating for $\GT$, for every
$\xi \in \Gamma_{T}^{\infty}(H_{\R})$ there exists a unique operator
$W(\xi) \in \GT$ such that
$$W(\xi)\Omega=\xi$$
W is called Wick product.

The right creation operator, $a^{*}_{r}(f)$, is defined by the following formula :
$$ a^{*}_{r}(f)(f_{1}\otimes \cdots \otimes f_{n})=f_{1}\otimes \cdots \otimes f_{n}\otimes f$$
We will also denote by $a_{r}(f)$ the right annihilation operator, which is its adjoint with respect to the T-scalar product, by $G_{r}(f)$ the right gaussian operator, and by $\Gamma_{T,r}(H_{\R})$ the von Neumann algebra generated by all right gaussians. It is easy to see that $\Gamma_{T,r}(H_{\R}) \subset \Gamma_{T}(H_{\R})^{\prime}$. Actually, by Tomita's theory, we have
$$\Gamma_{T,r}(H_{\R})=S\Gamma_{T}(H_{\R})S=\Gamma_{T}(H_{\R})^{\prime}$$
where $S$ is the anti linear operator on $\FT$ (which is actually an anti unitary) defined by

$$S(f_{1}\otimes \cdots \otimes f_{n})=f_{n}\otimes \cdots \otimes f_{1}$$
for any $f_{1}, \cdots,  f_{n}$ $\in H_{\R}$. Since $\Omega$ is also separating for $\Gamma_{T,r}(H_{\R})$ we can define the right Wick product, that will be denoted by $W_r (\xi)$. For any $\xi \in \Gamma_{T}^{\infty}(H_{\R})$ we have
$$(W(\xi))^* =W(S\xi) {\hskip 1cm}{\rm and}{\hskip 1cm} SW(\xi)S=W_{r}(S\xi)$$

Some particular cases of deformation have been
studied in the literature. Let $(q_{ij})_{i,\;j\in I}$ be a hermitian
matrix such that $\sup_{i,j}|q_{ij}|<1$. Define
$$Te_{i}\otimes e_{j}=q_{ij}e_{j}\otimes e_{i}$$
Then $T$ is a strictly contractive Yang-Baxter operator, and it is
tracial if and only if the $q_{ij}$  are real. Our deformed Fock space
is then a realisation of the following $q_{ij}$-relations :
$$a(e_{i})a^{*}(e_{j})-q_{ij}a^{*}(e_{j})a(e_{i})=\delta_{ij}$$
In the special case where all $q_{ij}$ are equal, we obtain the well
known q-relations.

Let us define the following selfadjoint unitary on the free Fock space :
$$\forall\;f_{1},\dots,f_{n} \in H_{\C}, \;\;\;U(f_{1}\otimes\dots\otimes f_{n})=f_{n}\otimes\dots\otimes f_{1}$$
Since $UP_{T}^{(n)}=P_{T}^{(n)}U$ (cf. \cite{K}), $U$ is also a
selfadjoint unitary on each T-Fock space.

Given vectors $f_{1},\dots,f_{n}$ in $H_{\R}$ we define :
$$a^{*}(f_{1}\otimes \dots \otimes f_{n})=a^{*}(f_{1})\dots
a^{*}(f_{n})\;\;\;{\rm and}\;\;\;a(f_{1}\otimes \dots \otimes
f_{n})=a(f_{1})\dots a(f_{n})$$
For $0 \le k \le n$, let $R_{n,k}$ be the operator on $H^{\otimes
  n}_{\C}$ given by
$$R_{n,k}=\sum_{\sigma \in S_{n}/S_{n-k}\times S_{k}}\varphi(\sigma^{-1})$$
where the sum runs over the representatives of the right cosets of
$S_{n-k}\times S_{k}$ in $S_{n}$ with minimal number of inversions.
Then
\begin{equation}P_{T}^{(n)}=R_{n,k}\left(P_{T}^{(n-k)}\otimes
  P_{T}^{(k)}\right
)\;\;\;{\rm and }\;\;\;\|R_{n,k}\|\le C_{q}
\end{equation}
where $C_{q}=\prod\limits_{n=1}^{\infty}(1-q^{n})^{-1}$ (cf. \cite{B2} and \cite{K}).
It follows that
\begin{equation}P_{T}^{(n)}\le C_{q} P_{T}^{(n-k)}\otimes P_{T}^{(k)} \end{equation}
It also follows that $a^{*}$, respectively $a$, extend linearly, respectively
antilinearly, and continuously to $H_{\C}^{\otimes n}$ for every $n \ge 1$. Then for each vector $\xi \in
H_{\C}^{\otimes n}$ we have
\begin{equation}\label{a}
\|a^{*}(\xi)\| \le C_{q}^{\frac{1}{2}}\|\xi\|_{T}\;\;\;{\rm and}\;\;\; (a^{*}(\xi))^{*}=a(U\xi).
\end{equation}
Let $n \ge 1$ and $1\le k \le n$, $H^{\otimes n-k}_{\C}\otimes H^{\otimes k}_{\C}$ will be the Hilbert tensor product of the Hilbert spaces $H^{\otimes k}_{\C}$ and $H^{\otimes n-k}_{\C}$ where both $H^{\otimes k}_{\C}$ and $H^{\otimes n-k}_{\C}$ are equipped with the T-scalar product.

\begin{lem}There is a positive constant $D_{q,n,k}$ such that
$$P_{T}^{(n-k)}\otimes P_{T}^{(k)} \le D_{q,n,k}P_{T}^{(n)}$$
Consequently for every $n \ge 1$ and $1 \le k \le n$, $H^{\otimes n}_{\C}$ and $H^{\otimes k}_{\C}\otimes H^{\otimes n-k}_{\C}$ are algebraically the same and their norms are equivalent.
\end{lem}

\noindent{ \bf Remark :} It is still not known whether one can choose
$D_{q,n,k}$ independent of $n$ and $k$.

\medskip

\proof It was shown in \cite{B1} that there is a positive constant $\omega(q)$ such that
$$P_{T}^{(n-1)}\otimes I \le \omega(q)^{-1}P^{(n)}_{T}$$
Since $U(P_{T}^{(n-1)}\otimes I)U=I \otimes P_{T}^{(n-1)}$ we also have
\begin{equation}I \otimes P_{T}^{(n-1)} \le \omega(q)^{-1}P^{(n)}_{T}\end{equation}
Fix some $k$, $2 \le k \le n-1$, using (3) and (4) we get :
\begin{eqnarray*}
P_{T}^{(n-k+1)}\otimes P_{T}^{(k-1)} & \le & C_{q}P_{T}^{(n-k)}\otimes I \otimes P_{T}^{(k-1)}\\
& \le & C_{q}\omega(q)^{-1}P_{T}^{(n-k)}\otimes P_{T}^{(k)}\\
\end{eqnarray*}
Thus by iteration it follows that for $0 \le k \le n$ :
\begin{equation} P_{T}^{(n-k)}\otimes P_{T}^{(k)} \le \omega(q)^{-1}(C_{q}\omega(q)^{-1})^{n-k}P_{T}^{(n)}\end{equation}
Since  $U(P_{T}^{(n-k)}\otimes P_{T}^{(k)})U=P_{T}^{(k)} \otimes P_{T}^{(n-k)}$ it follows from (6) that
\begin{equation*}  P_{T}^{(k)} \otimes P_{T}^{(n-k)} \le \omega(q)^{-1}(C_{q}\omega(q)^{-1})^{n-k}P_{T}^{(n)}\end{equation*}
Combining this last inequality and (6) we finally obtain :
\begin{equation} P_{T}^{(n-k)}\otimes P_{T}^{(k)} \le \omega(q)^{-1}(C_{q}\omega(q)^{-1})^{\min(k,n-k)}P_{T}^{(n)}\end{equation}
Then the desired result follows from (3) and (7).
\Eproof

\bigskip

For $k \ge 0$ let us now define on the family of finite linear
combinations of elementary tensors of length not less than $k$ the following operator $U_{k}$:
$$U_{k}(f_{1}\otimes \dots \otimes f_{n})=a^{*}(f_{1}\otimes \dots
\otimes f_{n-k})a(\overline{f}_{n-k+1}\otimes \dots \otimes
\overline{f}_{n})$$
where $\overline{\xi+i\eta}=\xi-i\eta$ for all $\xi,\; \eta \in
H_{\R}$.

\bigskip

Fix $n$ and $k$ with $n \ge k$. Let $\mathcal{J}:\;H_{\C}^{\otimes k}\rightarrow \overline{H_{\C}^{\otimes k}}$ be the conjugation (which is an anti isometry). For any $f_{1}, \cdots , f_{n}$, $\mathcal{J}$ is defined by $\mathcal{J}(f_{1}\otimes \cdots \otimes f_{n})=\overline{f_{1}}\otimes \cdots \otimes\overline{f_{n}}$.
It is clear that $U_{k}$ extends boundedly to $H_{\C}^{\otimes n-k}\otimes H_{\C}^{\otimes k}$ by the formula :
$$U_{k}=M(a^{*}\otimes a\mathcal{J})$$
where $M$ is the multiplication operator from $B(\FT)\otimes_{\min}B(\FT)$ to $B(\FT)$ defined by $M(A\otimes B)=AB$. Moreover, by (\ref{a}) we have
$$\|U_{k}\| \le \|M\|.\|a^{*}\otimes a\mathcal{J}\| \le C_{q}$$
where $U_{k}$ is viewed as an operator from $H_{\C}^{\otimes n-k}\otimes H_{\C}^{\otimes k}$ to $B(\FT)$.

\medskip

In the following lemma we state an extension of the Wick formula
(Theorem 3 in \cite{K}). We deduce it as an easy consequence of the
original Wick formula and of our previous discussion.

\begin{lem}Let $n \ge 1$ and $\xi \in \HC^{\otimes n}$, then $\HC^{\otimes n}\subset \Gamma_{T}^{\infty}(H_{\R})$ and we have the following Wick formula :
\begin{equation}\label{eqn}
W(\xi)=\sum_{k=0}^{n}U_{k}R_{n, k}^{*}(\xi)
\end{equation}
Moreover
\begin{equation}\label{boze}
\|\xi\|_{q} \le \|W(\xi)\| \le C_{q}^{\frac{3}{2}} (n+1)\|\xi\|_{q}
\end{equation}
\end{lem}

\noindent {\bf Remark} : (\ref{boze}) is the well known Bozejko's inequality discussed in \cite{B2} and \cite{K}, and which implies the ultracontractivity of the q-Ornstein Uhlenbeck semigroup. We include an elementary and simple proof.

\proof The usual Wick formula is the following (cf \cite{B2} and \cite{K}) : $\forall f_{1}, \dots , f_{n} \in \HC$ we have
$$W(f_{1}\otimes \dots \otimes f_{n})=\sum_{k=0}^{n}\sum_{\sigma \in S_{n}/S_{n-k}\times S_{k}}U_{k}\varphi(\sigma)(f_{1}\otimes \dots \otimes f_{n})$$
Hence (\ref{eqn}) holds for every $\xi \in \mathcal{A}_{n}=\{{\rm linear\; combinations\; of\; elementary\; tensors\; of\; length\; n}\}$. By Lemma 1 and our previous discussion, the right handside of (\ref{eqn}) is continuous from $H_{\C}^{\otimes n}$ to $B(\FT)$. Since $\Omega$ is separating, it follows that $H_{\C}^{\otimes n} \subset \Gamma_{T}^{\infty}(H_{\R})$ and that (\ref{eqn}) extends by density from $\mathcal{A}_{n}$ to $H_{\C}^{\otimes n}$. Actually, our argument shows that for any $\xi \in H^{\otimes n}$, $W(\xi)$ belongs to $C_{T}^{*}(H_{\R})$ which is the $C^{*}$-algebra generated by the T-gaussians.

\medskip

Since for any $\xi \in H_{\C}^{\otimes n}$, $W(\xi)\Omega =\xi$, the left inequality in (\ref{boze}) holds.
We have just showed that $W$ is bounded from $H_{\C}^{\otimes n}$ to $B(\FT)$. Hence, there is a constant $B_{q,n}$ such that for any $\xi \in H_{\C}^{\otimes n}$ we have $\|W(\xi)\| \le B_{q,n} \|\xi\|_{q}$. To end the proof of (\ref{boze}) we now give a precise estimate of $B_{q,n}$. Let $\xi \in H_{\C}^{\otimes n}$, by (\ref{eqn})  and (3) we have

\begin{equation}\label{est}
\|W(\xi)\| \le  \sum_{k=0}^{n}\|U_{k}R_{n,k}^{*}(\xi)\| \le
C_{q}\sum_{k=0}^{n}\|R_{n,k}^{*}(\xi)\|_{H_{\C}^{\otimes n-k}\otimes H_{\C}^{\otimes k}}
\end{equation}
It remains to compute the norm of $R_{n,k}^{*}$ as an operator from $H_{\C}^{\otimes n}$ to $H_{\C}^{\otimes n-k}\otimes H_{\C}^{\otimes k}$. Let $\eta \in H_{\C}^{\otimes n}$ we have, by (2) and (3)
\begin{eqnarray*}
\|R_{n,k}^{*}\eta\|^{2}_{H_{\C}^{\otimes n-k}\otimes H_{\C}^{\otimes k}}
& = & \Dual{P_{T}^{(n-k)}\otimes P_{T}^{(k)} R_{n,k}^{*}\eta}{R_{n,k}^{*}\eta}_{0}\\
& = & \Dual{P_{T}^{(n)}\eta}{R_{n,k}^{*}\eta}_{0}\le \|\eta\|_{T} \|R_{n,k}^{*}\eta\|_{T}
\end{eqnarray*}
On the other hand,

\begin{eqnarray*}
\|R_{n,k}^{*}\eta\|_{T}^{2} & = & \Dual{P_{T}^{(n)}R_{n,k}^{*}\eta}{R_{n,k}
^{*}\eta}_{0}\le C_{q} \Dual{P_{T}^{(n-k)}\otimes P_{T}^{(k)} R_{n,k}^{*}\eta}{R_{n,k}^{*}\eta}_{0}\\
& \le & C_{q}\Dual{P_{T}^{(n)}\eta}{R_{n,k}^{*}\eta}_{0}\\
& \le & C_{q}\|\eta\|_{T}\|R_{n,k}^{*}\eta\|_{T}
\end{eqnarray*}
Hence it follows that $\|R_{n,k}^{*}\eta\|_{T} \le C_{q}\|\eta\|_{T}$ and $\|R_{n,k}^{*}\eta\|^{2}_{H_{\C}^{\otimes n-k}\otimes H_{\C}^{\otimes k}} \le C_{q}\|\eta\|_{T}^{2}$.
Thus $\|R_{n,k}^{*}\|\le C_{q}^{\frac{1}{2}}$ as an operator from $H_{\C}^{\otimes n}$ to $H_{\C}^{\otimes n-k}\otimes H_{\C}^{\otimes k}$. From (\ref{est}) and this last estimate, follows the second inequality in (\ref{boze}).

\Eproof

\bigskip

The remainder of this section is devoted to a simple proof of the
non-injectivity of the free von Neumann algebra $\GO$ ($\dim H_{\R} \ge
2$). The main ingredient is the vector valued Bozejko inequality (Lemma
3 below), which is the free Fock space analogue of the corresponding
inequality for the free groups proved by Haagerup and Pisier in
\cite{HP} and extended by Buchholz in \cite{Bu2} (see also
\cite{Bu}). Note also that the inequality (11) below was first proved in
\cite{HP} in the case $n=1$ (i.e. for free gaussians) and that a
similar inequality holds for products of free gaussians (see
\cite{Bu2}).

\medskip

We will need the following notations : $(e_{i})_{i \in I}$ will denote
an orthonormal basis of $H_{\R}$, and for a multi-index $\underline{i}$ of
length $n$, $\underline{i}=(i_{1},\dots,i_{n})\in I^{n}$,
$e_{\underline{i}}=e_{i_{1}}\otimes \dots \otimes
e_{i_{n}}$. $(e_{\un{i}})_{|\un{i}|=n}$ is a real orthonormal basis of
$H^{\otimes n}_{\C}$ equipped with the free scalar product and
$(e_{\un{i}})_{|\un{i}|\ge 0}$
is a real orthonormal basis of the free Fock space.

\begin{lem}
Let $n\ge 1$, $K$ a complex Hilbert space and $(\alpha_{\un{i}})_{|\un{i}|=n}$ a finitely supported family of $B(K)$. Then :
\begin{equation}
\max_{0\le k \le n}\left\lbrace\left\|(\alpha_{\un{j},\;\un{l}})_{{|\un{j}|=n-k}\atop{|\un{l}|=k\;\;\;\;\;}}\right\|\right\rbrace
\le
\left\|\sum_{|\un{i}|=n}\alpha_{\un{i}}\otimes W(e_{\un{i}})\right\|
\le
(n+1)\max_{0\le k \le n}\left\lbrace\left\|(\alpha_{\un{j},\;\un{l}})_{{|\un{j}|=n-k}\atop{|\un{l}|=k\;\;\;\;\;}}\right\|\right\rbrace
\end{equation}
\end{lem}

\noindent{\bf Remark :} Since $(\alpha_{\un{i}})_{|\un{i}|=n}$ is
finitely supported the operator-coefficient matrix
$(\alpha_{\un{j},\;\un{l}})_{{|\un{j}|=n-k}\atop{|\un{l}|=k\;\;\;\;\;}}$
is a finite matrix, say a $r\times s$ matrix, and its norm is the
operator norm in $B(l_{2}^{s}(K),l_{2}^{r}(K))$.\\
\proof We write
$$\sum_{|\un{i}|=n}\alpha_{\un{i}}\otimes W(e_{\un{i}})=\sum_{k=0}^{n}F_{k}$$
where
$$F_{k}=\sum_{{|\un{j}|=n-k}\atop{|\un{l}|=k\;\;\;\;\;}}\alpha_{\un{j},\;\un{l}}\otimes a^{*}(e_{\un{j}})a(e_{\un{l}})$$
we have
$$F_{k}=(\dots I_{K}\otimes a^{*}(e_{\un{j}})
\dots)_{|\un{j}|=n-k}(\alpha_{\un{j},\;\un{l}}\otimes I_{\FO})_{{|\un{j}|=n-k}\atop{|\un{l}|=k\;\;\;\;\;}}
\left(
\begin{array}{c}
\vdots\\
I_{K}\otimes a(e_{\un{l}})\\
\vdots
\end{array}
\right)_{|\un{l}|=k}$$
that is, $F_{k}$ is a product of three matrices, the first is a row indexed by $\un{j}$, the third a column indexed by $\un{l}$. Note that

$$\|(\dots a^{*}(e_{\un{j}})\dots)_{|\un{j}|=n-k}\|^{2}=\|\sum_{|\un{j}|=n-k}a^{*}(e_{\un{j}})(a^{*}(e_{\un{j}}))^{*}\|=\|\sum_{|\un{j}|=n-k}a^{*}(e_{\un{j}})a(Ue_{\un{j}})\|$$
It is easy to see that $\sum\limits_{|\un{j}|=n-k}a^{*}(e_{\un{j}})a(Ue_{\un{j}})$ is the
orthogonal projection on $\bigoplus\limits_{p\ge n-k}H^{\otimes p}$.\\
Thus
$$\|(\dots a^{*}(e_{\un{j}})\dots)_{|\un{j}|=n-k}\|\le 1$$

Therefore

\begin{eqnarray*}
\|F_{k}\| & \le & \|(\dots I_{K}\otimes a^{*}(e_{\un{j}})
\dots)_{|\un{j}|=n-k}\|. \|(\alpha_{\un{j},\;\un{l}}\otimes I_{\FO})_{{|\un{j}|=n-k}\atop{|\un{l}|=k\;\;\;\;\;}}\|.
\left\|
\left(
\begin{array}{c}
\vdots\\
I_{K}\otimes a(e_{\un{l}})\\
\vdots
\end{array}
\right)_{|\un{l}|=k}
\right\|\\
& \le & \|(\dots a^{*}(e_{\un{j}})\dots)_{|\un{j}|=n-k}\|. \|(\alpha_{\un{j},\;\un{l}})_{{|\un{j}|=n-k}\atop{|\un{l}|=k\;\;\;\;\;}}\|.
\|(\dots a^{*}(Ue_{\un{l}})\dots)_{|\un{l}|=k}\|\\
& \le & \|(\alpha_{\un{j},\;\un{l}})_{{|\un{j}|=n-k}\atop{|\un{l}|=k\;\;\;\;\;}}\|
\end{eqnarray*}

It follows that
$$\|\sum_{|\un{i}|=n}\alpha_{\un{i}}\otimes W(e_{\un{i}})\| \le
\sum_{k=0}^{n}\|F_{k}\| \le (n+1)\max_{0\le k \le n}\|(\alpha_{\un{j},\;\un{l}})_{{|\un{j}|=n-k}\atop{|\un{l}|=k\;\;\;\;\;}}\|$$
\allowdisplaybreaks{
To prove the first inequality, fix $0\le k_{0} \le n$ and consider
$(v_{\un{p}})_{|\un{p}|=k_{0}}$ such that $\sum\limits_{|\un{p}|=k_{0}}\|v_{\un{p}}\|^{2}< +\infty$. Let
$\eta =\sum\limits_{|\un{p}|=k_{0}}v_{\un{p}}\otimes Ue_{\un{p}}$. We have :
\begin{eqnarray*}
\|\sum_{|\un{i}|=n}\alpha_{\un{i}}\otimes W(e_{\un{i}})\eta\|^{2} & =
& \sum_{k=0}^{n}\|F_{k}\eta\|^{2} \ge \|F_{k_{0}}\eta\|^{2}\\
& = & \|\sum_{{|\un{j}|=n-k_{0}}\atop{|\un{l}|=k_{0}\;\;\;\;\;}}\alpha_{\un{j},\;\un{l}}v_{\un{l}}\otimes e_{\un{j}}\|^{2}\\
& = & \sum_{|\un{j}|=n-k_{0}}\|\sum_{|\un{l}|=k_{0}}\alpha_{\un{j},\;\un{l}}v_{\un{l}}\|^{2}\\
& = &\left\|(\alpha_{\un{j},\;\un{l}})_{{|\un{j}|=n-k_{0}}\atop{|\un{l}|=k_{0}\;\;\;\;\;}}
\left(
\begin{array}{c}
\vdots\\
v_{\un{l}}\\
\vdots
\end{array}
\right)_{|\un{l}|=k_{0}}
\right\|^{2}\\
\end{eqnarray*}
Then the result follows.
\Eproof}

\bigskip

Using Lemma 3, it is now easy to prove that $\GO$ is not injective as
soon as $\dim H_{\R} \ge 2$. Suppose that $\GO$ is injective and
$\dim H_{\R} \ge 2$. Choose two orthonormal vectors
$e_{1}$ and $e_{2}$ in $H_{\R}$. For $n\ge 1$ we have by
semi-discreteness (which is equivalent to the injectivity):
{\allowdisplaybreaks \begin{equation*}
\tau\left(\sum_{|\un{i}|=n}W(e_{\un{i}})^{*}W(e_{\un{i}})\right)  \le
     \|\sum_{|\un{i}|=n}\overline{W(e_{\un{i}})}\otimes W(e_{\un{i}})\|
\end{equation*}
where in the above sums, the index $\un{i} \in \{1,2\}^{n}$. However,
\begin{eqnarray*}
\tau\left(\sum_{|\un{i}|=n}W(e_{\un{i}})^{*}W(e_{\un{i}})\right) & = &
\sum_{|\un{i}|=n}\Dual{W(e_{\un{i}})\Omega}{W(e_{\un{i}})\Omega}_{0}\\
 & = & \sum_{|\un{i}|=n}\|e_{\un{i}}\|^{2} = 2^n
\end{eqnarray*}
On the other hand, by Lemma 3,
\begin{eqnarray*}
\|\sum_{|\un{i}|=n}\overline{W(e_{\un{i}})}\otimes W(e_{\un{i}})\| &
\le & (n+1)\max_{0\le k \le
    n}\lbrace\|(\overline{W(e_{\un{j},\;\un{l}})})_{{|\un{j}|=n-k}\atop{|\un{l}|=k\;\;\;\;\;}}\|\rbrace\\
 & \le & (n+1)\left(\sum_{|\un{i}|=n}\|W(e_{\un{i}})\|^{2}\right)^{\frac{1}{2}}\\
 & \le & (n+1)(2^{n}(n+1)^{2})^{\frac{1}{2}}\\
 & \le & (n+1)^{2}2^{\frac{n}{2}}\\
\end{eqnarray*}
Combining the preceding inequalities, we get $2^n \le (n+1)^{2}2^{\frac{n}{2}}$
which yields a contradiction for sufficiently large $n$. Therefore, $\GO$ is not injective if $\dim H_{\R} \ge 2$.}

\section{Generalized Haagerup-Bo$\rm{ \bf \dot{z}}$ejko inequality and non injectivity of $\GT$}

In the following we state and prove  the generalized inequality (1). It actually solves a question of Marek Bozejko ( in \cite{B2} page 210)
whether it is possible to find an operator coefficient version
of the following inequality (this is inequality (\ref{boze}) in Lemma 2):
\begin{equation}\label{boz}
\|\sum_{|\un{i}|=n}\alpha_{\un{i}}e_{\un{i}}\| \le
\|\sum_{|\un{i}|=n}\alpha_{\un{i}}W(e_{\un{i}})\| \le
C_{q}^{\frac{3}{2}}(n+1)\|\sum_{|\un{i}|=n}\alpha_{\un{i}}e_{\un{i}}\|
\end{equation}
where $(\alpha_{\un{i}})_{\un{i}}$ is a finitely supported family of
 complex numbers. Inequality (12) was  proved in \cite{B2} for the
q-deformation, and generalized in \cite{K} for the Yang-Baxter deformation.

 First, we need to recall some basic notions from operator space theory. We refer to \cite{ER} and \cite{P} for more information.\\
Given $K$ a complex Hilbert space, we can equip $K$ with the column, respectively the row, operator space structure denoted by $K_{c}$, respectively $K_{r}$, and defined by
$$K_{c}=B(\C,K){\hskip 1cm}{\rm and}{\hskip 1 cm}K_{r}=B(K^{*},\C).$$
Moreover, we have $K_c^*=\overline K_r$ as operator spaces.

Given two operator spaces $E$ and $F$, let us briefly recall the definition of the Haagerup tensor product of $E$ and $F$. $E\otimes F$ will denote the algebraic tensor product of $E$ and $F$. For $n\ge 1$ and $x=(x_{i,j})$ belonging to $M_{n}(E\otimes F)$ we define
$$\|x\|_{(h,n)}=\inf\{\|y\|_{M_{n,r}(E)}\|z\|_{M_{r,n}(F)}\}$$
where the infimum runs over all $r\ge 1$ and all decompositions of $x$ of the form
$$x_{i,j}=\sum_{k=1}^{r}y_{i,k}\otimes z_{k,j}.$$
By Ruan's theorem, this sequence of norms define an operator space structure on the completion of $E \otimes F$ equipped with $\|\;.\;\|_{h}=\|\;.\;\|_{(h,1)}$. The resulting operator space, which is called the Haagerup tensor product of $E$ and $F$ is denoted by $E\otimes_{h} F$.\\
In this setting, a bilinear map $u:\;E\times F\rightarrow B(K)$ is said to be completely bounded, in short c.b, if and only if the associated linear map $\hat{u}:\; E\otimes F \rightarrow B(K)$ extends completely boundedly to $E\otimes_{h} F$. We define $\|u\|_{cb}=\|\hat{u}\|_{cb}$. This notion goes back to Christensen and Sinclair \cite{CS}.

 We will often use  the following classical identities for hilbertian
operator spaces :
$$ K_c \otimes_{\min} H_r = K_c \otimes_{h} H_r =
\mathcal K(\overline H, K),$$
where $\mathcal K$ stands for the compact operators and
$$  K_c \otimes_{\min} H_c = K_c \otimes_{h} H_c =(K \otimes_{2} H)_c$$
and similarly for rows using duality.

There is another notion of complete boundedness for bilinear maps, called jointly complete boundedness. Let $E$, $F$ be operator spaces, $K$ a complex Hilbert space, and $u:\;E\times F\rightarrow B(K)$ a bilinear map. $u$ is said to be jointly completely bounded (in short j.c.b) if and only if for any ${\rm C}^{*}$-algebras $B_{1}$ and $B_{2}$, $u$ can be boundedly extended to a bilinear map $(u)_{B_{1},B_{2}}:\;E\otimes_{\min}B_{1}\times F\otimes_{\min}B_{2} \rightarrow B(K)\otimes_{\min}B_{1}\otimes_{\min}B_{2}$ taking $(e\otimes b_{1},\; f\otimes b_{2})$ to $u(e,\;f)\otimes b_{1}\otimes b_{2}$. We put $\|u\|_{jcb}=\sup\limits_{B_{1},\;B_{2}}\|(u)_{B_{1},B_{2}}\|$. Observe that in this definition $B_{1}$ and $B_{2}$ can be replaced by operator spaces.

We will need the fact that every bilinear c.b map is a j.c.b map with $\|u\|_{jcb} \le \|u\|_{cb}$. Let $K$ be a complex Hilbert space and $u:\; B(K)\times K_{c} \rightarrow K_{c}$ the bilinear map taking $(\varphi,\;k)$ to $\varphi(k)$. Then it is easy to see that $u$ is a norm one bilinear cb map.\\
To simplify our notations, $H_{\C}$ will be, most of the time, replaced by $H$ in the rest of this section. For the same reason we will denote by $H_c ^{\otimes n}$ (respectively $H_r ^{\otimes n}$) the column Hilbert space $(H_{\C}^{\otimes n})_{c}$ (respectively the row Hilbert space $(H_{\C}^{\otimes n})_{r}$).

\begin{lem}
Let $n\ge 1$. The mappings  $a^{*}:\;H_{c}^{\otimes n} \rightarrow B(\FT)$ and $a : \overline H_r^{\otimes n} \to  B(\FT)$ are completely bounded with cb-norms less than $\sqrt C_q$.
\end{lem}

\proof
Let us start with the proof of the statement concerning $a^{*}$. Let $n\ge 1$, $K$ a complex Hilbert space and $(\alpha_{\un{i}})_{|\un{i}|=n}$ a finitely supported family of $B(K)$ such that
$$\|\sum_{|\un{i}|=n}\alpha_{\un{i}}\otimes e_{\un{i}}\|_{B(K)\otimes_{\rm min}H_{c}} <1.$$

Then, since the maps $a^*(e_{\un{i}})$ acts diagonally with respect to degrees
of tensors in $\FT$,
$$\|\sum_{|\un{i}|=n}\alpha_{\un{i}}\otimes a^{*}(e_{\un{i}})\|_{B(K)\otimes_{\rm min}B(\FT)}=\sup_{k\ge 0}\|\sum_{|\un{i}|=n}\alpha_{\un{i}}\otimes a^{*}(e_{\un{i}})\|_{B(K)\otimes_{\rm min}B(H^{\otimes k},H^{\otimes n+k})}$$

To compute the right term, fix $k\ge 0$ and let
$(\xi_{\un{j}})_{|\un{j}|=k}$ be
a finitely supported family of vectors in $K$ such that
$$\|\sum_{|\un{j}|=k}\xi_{\un{j}}\otimes
e_{\un{j}}\|_{K\otimes_{2}H^{\otimes k}}<1.$$

By (3) we have
$$\|\sum_{\un{i},\;\un{j}}\alpha_{\un{i}}(\xi_{\un{j}})\otimes e_{\un{i}} \otimes e_{\un{j}}\|_{K\otimes_{2}H^{\otimes n+k}} \le C_{q}^{\frac{1}{2}}\|\sum_{\un{i},\;\un{j}}\alpha_{\un{i}}(\xi_{\un{j}})\otimes e_{\un{i}} \otimes e_{\un{j}}\|_{K\otimes_{2}H^{\otimes n}\otimes_{2}H^{\otimes k}}.$$

Let $u:\; B(K)\times K_{c} \rightarrow K_{c}$ given by $(\varphi,\xi) \mapsto \varphi(\xi)$. Recall that  $\|u\|_{cb}= 1$. Consequently, $\|u\|_{jcb}\le 1$. Therefore, we deduce

\begin{eqnarray*}
\|\sum_{\un{i},\;\un{j}}\alpha_{\un{i}}(\xi_{\un{j}})\otimes e_{\un{i}}
\otimes e_{\un{j}}\|_{K\otimes_{2}H^{\otimes n}\otimes_{2}H^{\otimes k}}
& = & \|\sum_{\un{i},\;\un{j}}\alpha_{\un{i}}(\xi_{\un{j}})\otimes e_{\un{i}}
\otimes e_{\un{j}}\|_{K_{c}\otimes_{\min}H^{\otimes n}_{c}
\otimes_{\min}H^{\otimes k}_{c}}\\
& = & \|(u)_{H_{c}^{\otimes n},\;H_{c}^{\otimes k}}(\sum_{\un{i}}\alpha_{\un{i}}\otimes e_{\un{i}},\;\sum_{\un{j}}\xi_{\un{j}} \otimes e_{\un{j}})\|\\
& \le & \|u\|_{jcb}\|\sum_{\un{i}}\alpha_{\un{i}}\otimes e_{\un{i}}\|_{B(K)\otimes_{\rm min}H_{c}^{\otimes n}}\|\sum_{\un{j}}\xi_{\un{j}} \otimes e_{\un{j}}\|_{K_{c}\otimes_{\rm min}H_{c}^{\otimes k}}\\
& \le & 1
\end{eqnarray*}

\medskip

By the result just proved, for any complex Hilbert space $K$ and for any finitely supported family $(\alpha_{\un{i}})_{|\un{i}|=n}$ of $B(K)$ we have
$$\|\sum_{|\un{i}|=n}\alpha_{\un{i}}\otimes a^* (e_{\un{i}})\|_{B(K)\otimes_{\min} B(\FT)}
\le \sqrt{C_{q}} \|\sum_{|\un{i}|=n}\alpha_{\un{i}}\otimes
e_{\un{i}}\|_{B(K)\otimes_{\min} H_{c} ^{\otimes n}}$$
Taking adjoints on both sides we get
$$\|\sum_{|\un{i}|=n}\alpha_{\un{i}}^{*}\otimes a(Ue_{\un{i}})\|_{B(K)\otimes_{\min} B(\FT)}
\le \sqrt{C_{q}} \|\sum_{|\un{i}|=n}\alpha_{\un{i}}^{*}\otimes
\overline{e_{\un{i}}}\|_{B(K)\otimes_{\min} \overline{H}_{r} ^{\otimes n}}$$

Changing $\alpha_{\un{i}}^{*}$ to $\alpha_{\un{i}}$ and using the fact that $U$
(reversing the order of tensor) is a complete isometry on  $H_r^{\otimes n}$,
 we get that for any finitely supported family $(\alpha_{\un{i}})_{|\un{i}|=n}$ of $B(K)$ we have

$$\|\sum_{|\un{i}|=n}\alpha_{\un{i}}\otimes a(\overline{e_{\un{i}}})\|_{B(K)\otimes_{\min} B(\FT)}
\le \sqrt{C_{q}} \|\sum_{|\un{i}|=n}\alpha_{\un{i}}\otimes
\overline{e_{\un{i}}}\|_{B(K)\otimes_{\min} \overline{H}_{r} ^{\otimes n}}.$$
In other words,
$$ a : \overline H_r^{\otimes n} \to  B(\FT)$$
is also completely bounded with norm less than $\sqrt C_q$.
\Eproof

\begin{cor}\label{estuk} For any $n \ge 0$, and any $k\in\{0...n\}$,
$$U_{k}:\; H^{\otimes n-k}_{c}\otimes_{h} H^{\otimes k}_{r}
\rightarrow B(\FT)$$
is  completely bounded with cb-norm less  than $C_{q}$.
\end{cor}

\proof Let us denote by $M$ the
multiplication map $B(\FT)\otimes_h B(\FT)\rightarrow B(\FT)$ given by $A\otimes B \mapsto AB$, $M$ is obviously
completely contractive. We have the formula
$$U_k=M(a^*\otimes a\mathcal{J})$$
if $\mathcal{J}: H^{\otimes k}\to \overline H^{\otimes k}$ is the conjugation
(which is a complete isometry).
By injectivity of the Haagerup tensor product and by Lemma 4 we deduce that
$$\|a^{*}\otimes a \mathcal{J}\|_{cb}\le  C_q$$
 Then
$$\|U_k\|_{cb} \le \|M\|_{cb} \|a^{*}\otimes a \mathcal{J}\|_{cb}\le  C_q$$

\Eproof

Recall that, by definition, $\Gamma^{\infty}_{T}(H_{\R})$ is identified with $\GT$ by the mapping sending $\xi$ to $W(\xi)$. Thus $\Gamma^{\infty}_{T}(H_{\R})$ inherits the operator space structure of $\GT$. In particular for all $n \ge 0$, $H^{\otimes n}$ will be equipped with the operator space structure of $E_{n}=\{W(\xi),\;\xi\in H^{\otimes n}\}$.

Theorem 1 below was first obtained via elementary, but long, computations. In the version presented here, we have chosen to follow an approach indicated to us by Eric Ricard. This approach is much more transparent but involves some notions of operator space theory.

\bigskip

\begin{thm}
Let $K$ be a complex Hilbert space. Then for all $n\ge 0$ and for all $\xi \in B(K)\otimes_{\min} H^{\otimes n}$ we have
\begin{equation}
\max_{0\le k \le n}\|( Id\otimes R_{n,\;k}^{*})(\xi)\|
\le \|({ Id}\otimes W)(\xi)\|_{\min}\le C_{q}(n+1)
\max_{0\le k \le n}\|({ Id}\otimes R_{n,\;k}^{*})(\xi)\|
\end{equation}
where ${ Id}$ denotes the identity mapping of $B(K)$, and where the
 norm  $\|({ Id}\otimes R_{n,\;k}^{*})(\xi)\|$ is that of $B(K)\otimes_{\min}H^{\otimes n-k}_{c}\otimes_{\min}H^{\otimes k}_{r}$.
\end{thm}

\proof
For the second  inequality, we use the Wick formula :
$$W\left|_{H^{\otimes n}} \right. =\sum_{k=0}^n U_k \, R_{n,k}^*.$$
Let $\xi \in B(K)\otimes_{\min} H^{\otimes n}$, then by corollary \ref{estuk}
$$\|({ Id}\otimes W)(\xi)\|_{\min}\le C_{q}
\sum_{k=0}^n\|({ Id}\otimes R_{n,\;k}^{*})(\xi)\|$$
which yields the majoration.

\smallskip

 For the minoration, for $x\in H_{c}^{\otimes n-k}\otimes
H_{r}^{\otimes k}\subset B(\overline H^{\otimes k},H^{\otimes n-k})$, we claim that
\begin{equation}\label{coucou}
P_{n-k} U_k (x)\left|_{H^{\otimes^k}}\right. = x(U\mathcal{J})
\end{equation}
where $P_{n-k}$ is the projection on tensors of rank $n-k$ in $\FT$.
Assuming this claim and recalling that $U$ and $\mathcal{J}$ are (anti)-isometry,
we get that for any $x\in B(K)
\otimes_{\min} H_c^{\otimes n-k}\otimes_{\min} H_r^{\otimes k}$
$$\|x\|_{B(K)
\otimes_{\min} H_c^{\otimes n-k}\otimes_{\min} H_r^{\otimes k}}\le
\|P_{n-k}\|_{B(\FT)} \|(Id\otimes U_k)(x)\|_{B(K)\otimes_{\min}B(\FT)}$$
The conclusion follows applying this inequality to
$x=({ Id}\otimes R_{n,\;k}^{*})(\xi)$

To prove (\ref{coucou}), it suffices to consider an elementary
tensor product with entries in any basis of $H$, say
$x=e_{\un{i}}\otimes e_{\un{j}}$. Consider $e_{\un{l}}\in H^{\otimes
k}$, a length argument gives that $a(\mathcal{J}e_{\un j}).e_{\un
l}$ is of the form $\lambda \Omega$, with
$$\lambda= \Dual{a(\mathcal{J}e_{\un j}).e_{\un l}} \Omega = \Dual{e_{\un l}}{\mathcal{J}U
e_{\un j}}$$
We deduce that
$$P_{n-k} U_k (e_{\un{i}}\otimes
e_{\un{j}}) . e_{\un l}= \Dual{e_{\un l}}{U\mathcal{J}
e_{\un j}} \,e_{\un{i}}.$$
On the other hand, viewing $x$ as an operator, we compute
$$ x(\mathcal{J}U).e_{\un l}= x.(\mathcal{J}Ue_{\un l})=\Dual{e_{\un j}}{\mathcal{J}Ue_{\un l}}
 \, e_{\un{i}}$$
But since $U$ is unitary and $\mathcal{J}$ antiunitary,
$$\Dual{e_{\un j}}{\mathcal{J}Ue_{\un l}}= \Dual{e_{\un l}}{U\mathcal{J}e_{\un j}}$$
This ends the proof.
\Eproof

\bigskip

The following theorem is our main result.

\begin{thm}
$\GT$ is not injective as soon as  $\dim(H_{\R})\ge 2$.
\end{thm}

\proof Let $d\le\dim H_{\R} $.
For all $n\ge0$, $(\xi_{\un{i}})_{|\un{i}|=n}$ will denote a real orthonormal
family of $H^{\otimes n}$  equipped with the T-scalar product of cardinal
$d^n$. For example one can take
$\xi_{\un{i}}=(P_{T}^{(n)})^{-\frac{1}{2}}e_{\un{i}}$.

Suppose that $\GT$ is injective.
Fix $n\ge 1$. By injectivity we have,
\begin{equation*}
\tau(\sum_{|\un{i}|=n}W(\xi_{\un{i}})^{*}W(\xi_{\un{i}}))  \le
\|\sum_{|\un{i}|=n}\overline{W(\xi_{\un{i}})}\otimes W(\xi_{\un{i}})\|
\end{equation*}

It is clear that
$$\tau(\sum_{|\un{i}|=n}W(\xi_{\un{i}})^{*}W(\xi_{\un{i}}))= d^n$$
On the other hand, applying  twice (13) consecutively

\begin{equation*}
\|\sum_{|\un{i}|=n}\overline{W(\xi_{\un{i}})}\otimes
W(\xi_{\un{i}})\|
 \le  (n+1)^2C_{q}^2\max_{0\le k,k' \le n}\lbrace
\|\sum_{|\un i|=n}
\overline{R_{n,k'}^*(\xi_{\un i})}\otimes R_{n,k}^*(\xi_{\un i})\|\rbrace
\end{equation*}
The norms are computed in $ \overline H_c^{\otimes n-k'}\otimes_{\min}
\overline H_r^{\otimes k'}
\otimes_{\min}H_c^{\otimes n-k}\otimes_{\min} H_r^{\otimes k}$ for fixed $k$
and $k'$. We can rearrange this tensor product
and use the
comparison with the Hilbert Schmidt norm :
Let $t=\sum_{|\un i|=n}
\overline{R_{n,k'}^*(\xi_{\un i})}\otimes R_{n,k}^*(\xi_{\un i})$,
 \begin{eqnarray*}
\|t\|_{\overline H_c^{\otimes n-k'}\otimes_{\min}
\overline H_r^{\otimes k'}
\otimes_{\min}H_c^{\otimes n-k}\otimes_{\min} H_r^{\otimes k}}
 &=&\|t\|_{
(\overline H^{\otimes n-k'} \otimes_2 H^{\otimes n-k})_c\otimes_{\min}
(\overline H^{\otimes k'} \otimes_2 H^{\otimes k})_r}\\
&\le& \|t\|_{
(\overline H^{\otimes n-k'} \otimes_2 H^{\otimes n-k})\otimes_{2}
(\overline H^{\otimes k'} \otimes_2 H^{\otimes k})} \\
&\le &  \|t\|_
{\overline H^{\otimes n-k'}\otimes_{2}
\overline H^{\otimes k'}
\otimes_{2}H^{\otimes n-k}\otimes_{2} H^{\otimes k}}
\end{eqnarray*}
Finally, we use the estimates on $R_{n,k}^*$ :
\begin{eqnarray*}
\|t\|_{\overline H_c^{\otimes n-k'}\otimes_{\min}
\overline H_r^{\otimes k'}
\otimes_{\min}H_c^{\otimes n-k}\otimes_{\min} H_r^{\otimes k}}
 &\le &  \|t\|_
{\overline H^{\otimes n-k'}\otimes_{2}
\overline H^{\otimes k'}
\otimes_{2}H^{\otimes n-k}\otimes_{2} H^{\otimes k}}\\
& \le & C_q \|\sum_{|\un i|=n}
\overline{\xi_{\un i}}\otimes \xi_{\un i}\|_{\overline H^n\otimes_2 H^n}
\end{eqnarray*}
But by the choice of $\xi_{\un i}$ : $\|\sum_{|\un i|=n}
\overline{\xi_{\un i}}\otimes \xi_{\un i}\|_{\overline H^n\otimes_2 H^n}=
d^{n/2}$.

Combining all inequalities above, we deduce
$$d^{n}  \le  C_{q}^{3}(n+1)^{2}d^{n/2}$$
which yields a contradiction when $n$ tends to infinity as soon as $d\ge 2$.
\Eproof

\medskip

Let $C_{T}^{*}(H_{\R})$ be the $C^{*}$-algebra
generated by all gaussians $G(f)$ for $f\in H_{\R}$. The preceding theorem
implies directly that $C_{T}^{*}(H_{\R})$ is not nuclear as soon as
$\dim(H_{\R})\ge 2$ (cf. \cite{CE} Corollary 6.5). Actually the preceding argument can be modified to prove that $C_{T}^{*}(H_{\R})$ does not have the weak expectation property as soon as $\dim H_{\R} \ge 2$. Recall that a $C^{*}$-algebra $A$ has the weak expectation property
(WEP in short)
if and only if the canonical inclusion $A \rightarrow A^{**}$
factorizes completely contractively through $B(K)$
for some complex Hilbert space $K$. By the results of Haagerup
(cf. \cite{P} Chapter 15) a $C^{*}$-algebra $A$ has the WEP if and only if
for all finite family $x_{1},\hdots, x_{n}$ in $A$
\begin{equation}\label{Hag}
\|\sum_{i=1}^{n} x_{i}\otimes \overline{x_{i}}\|_{A\otimes_{\max}\overline{A}}=\|\sum_{i=1}^{n} x_{i}\otimes \overline{x_{i}}\|_{A\otimes_{\min}\overline{A}}
\end{equation}

\begin{cor} $C_{T}^{*}(H_{\R})$ does not have the WEP as soon as $\dim H_{\R} \ge 2$.
\end{cor}

\proof Let us use the same notations as in the preceding proof and suppose that $C_{T}^{*}(H_{\R})$ has the WEP. Fix $n \ge 1$, by (\ref{Hag}) we have

\begin{equation}\label{WEP}
\|\sum_{|\un{i}|=n} W(\xi_{\un{i}})\otimes \overline{W(\xi_{\un{i}})}\|_{C_{T}^{*}(H_{\R})\otimes_{\max}\overline{C_{T}^{*}(H_{\R})}}\le
\|\sum_{|\un{i}|=n} W(\xi_{\un{i}})\otimes \overline{W(\xi_{\un{i}})}\|_{C_{T}^{*}(H_{\R})\otimes_{\min}\overline{C_{T}^{*}(H_{\R})}}
\end{equation}
To estimate from below the left handside of (\ref{WEP}) observe that $\Phi :\;\overline{C_{T}^{*}(H_{\R})} \rightarrow C_{T}^{*}(H_{\R})^{\prime}$ taking $\overline{W(\xi)}$ to $\mathcal{J}UW(\xi)\mathcal{J}U=W_{r}(\mathcal{J}U\xi)$ is a $*$- representation. Thus

\begin{eqnarray*}
\|\sum_{|\un{i}|=n} W(\xi_{\un{i}})\otimes \overline{W(\xi_{\un{i}})}\|_{C_{T}^{*}(H_{\R})\otimes_{\max}\overline{C_{T}^{*}(H_{\R})}} & = & \|\sum_{|\un{i}|=n} W(\xi_{\un{i}})\otimes W_{r}(\mathcal{J}U\xi_{\un{i}})\|_{C_{T}^{*}(H_{\R})\otimes_{\max}C_{T}^{*}(H_{\R})^{\prime}}\\
& \ge & \|\sum_{|\un{i}|=n} W(\xi_{\un{i}})W_{r}(\mathcal{J}U\xi_{\un{i}})\|_{B(\FT)}\\
& \ge & \sum_{|\un{i}|=n}\Dual{\mathcal{J}U\xi_{\un{i}}}{W(\xi_{\un{i}})^{*}\Omega}_{T}\\
& \ge & \sum_{|\un{i}|=n}\Dual{\mathcal{J}U\xi_{\un{i}}}{W(\mathcal{J}U\xi_{\un{i}})\Omega}_{T}\\
& \ge & \sum_{|\un{i}|=n}\|\mathcal{J}U\xi_{\un{i}}\|^{2}_{T}=d^n
\end{eqnarray*}
Then we can finish the proof as for  Theorem 2.
\Eproof

\noindent{\bf Remark :}
 Non nuclearity of $C_{T}^{*}(H_{\R})$ is equivalent to the fact that
$C_{T}^{*}(H_{\R})$ does
not have the completely positive approximation
property as soon as $\dim(H_{\R})\ge
2$. However it is possible to prove that $C_{T}^{*}(H_{\R})$ has
the metric approximation property, by truncation of the Ornstein-Uhlenbeck
semigroup. Arguing by duality and interpolation, it is not difficult
to show that $L^{p}\left(\GT\right)$ has the metric approximation
property for $1\le p< \infty$. However, at the time of this writing, we are not able to prove that $C_{T}^{*}(H_{\R})$ has the completely bounded approximation property.

\section{The case of the $q-$Araki-Woods algebras}

For this last section we mainly refer to \cite{H} where the
 $q-$Araki-Woods algebras are defined as a generalization of the $q-$
 deformed case of Bo\.zejko and Speicher on the one hand, and the quasi-free case of
 Shlyakhtenko (cf. \cite{S}) on the other. More precisely, let $H_{\R}$ be a real
 Hilbert space given with $U_{t}$, a strongly continuous group of orthogonal
 transformations on $H_{\R}$. $U_{t}$ can be extended to a unitary group
 on the complexification $H_{\C}$. Let $A$ be its positive
 non-singular generator on $H_{\C}$ : $U_{t}=A^{it}$. A new scalar
 product $\Dual{\;.\;}{\;.\;}_{U}$ is defined on $H_{\C}$ by the
 following relation :
$$\Dual{\xi}{\eta}_{U}=\Dual{2A(1+A)^{-1}\xi}{\eta}$$
We will denote by $H$ the completion of $H_{\C}$ with
respect to this new scalar product.

\medskip

For a fixed $q\in ]-1,1[$, we now consider the $q-$deformed Fock space
associated with $H$ and we denote it by $\FQH$. Recall that it is
the Fock space with the following Yang-Baxter deformation $T$ defined
by :
\fonction{T}{H\otimes H}{H\otimes H}{\xi\otimes \eta}{q\eta\otimes \xi}
Or equivalently, for every $n\ge 2$ and $\sigma \in S_{n}$ we have
$$\varphi(\sigma)=q^{i(\sigma)}U_{\sigma}$$
where $i(\sigma)$ denotes the number of inversions of the permutation
$\sigma$ and $U_{\sigma}$ is the unitary on $H^{\otimes n}$ defined by
$$U_{\sigma}(f_{1}\otimes \dots \otimes f_{n})=f_{\sigma^{-1}(1)}\otimes
\dots \otimes f_{\sigma^{-1}(n)}$$

\medskip

\noindent In this setting, the $q-$Araki-Woods algebra is the following von
Neumann algebra
$$\GQt=\{G(h), h\in H_{\R}\}^{\prime\prime}\subset B\left(\FQ\right)$$
Let $H_{\R}^{\prime}=\{g \in H, \Dual{g}{h}_{U} \in \R {\rm\; for\; all\;} h \in
H_{\R}\}$ and
$$\GQtr=\{G_{r}(h), h\in H_{\R}^{\prime}\}^{\prime\prime}$$
where $G_{r}(h)$ is the right gaussian corresponding to the right
creation operator.\\
Since $\GQtr \subset \GQt^{\prime}$,
$\overline{H_{\R}+iH_{\R}}=H$ and
$\overline{H_{\R}^{\prime}+iH_{\R}^{\prime}}=H$ (cf. \cite{S}), it is
easy to deduce that $\Omega$ is cyclic and separating for both $\GQt$ and
$\GQtr$. So
Tomita's theory can apply : recall that the anti-linear operator $S$ is
the closure of the operator defined by :
$$S(x\Omega)=x^{*}\Omega\;\;\;{\rm\; for\; all\;} x\in \GQt$$
Let $S=J\Delta^{\frac{1}{2}}$ be its polar decomposition. $J$ and
$\Delta$ are called respectively the modular conjugation and the
modular operator. The following explicit formulas hold  (cf. \cite{H}
and \cite{S})
$$S(h_{1}\otimes \dots \otimes h_{n})=h_{n}\otimes \dots \otimes
h_{1}\;\;\;{\rm\; for\; all\;} h_{1},\dots,h_{n} \in H_{\R}$$
$\Delta$ is the closure of the operator
$\bigoplus\limits_{n=0}^{\infty}(A^{-1})^{\otimes n}$ and
$$J(h_{1}\otimes \dots \otimes h_{n})=A^{-\frac{1}{2}}h_{n}\otimes
\dots \otimes A^{-\frac{1}{2}}h_{1}\;\;\;{\rm\; for\; all\;}
h_{1},\dots,h_{n} \in H_{\R}\cap {\rm dom}A^{-\frac{1}{2}}$$
By Tomita's theory, we have
$$\GQt^{\prime}=J\GQt J$$
Let $h\in H_{\R}$, as in \cite{S} we have $Jh \in H_{\R^{\prime}}$,
then, since $\Omega$ is separating for $\GQtr$, we obtain that
$JG(h)J=G_{r}(Jh) \in \GQtr$, so that
$$\GQt^{\prime}=\GQtr$$
Moreover, if $\xi \in \GQt\Omega$, then $J\xi \in \GQtr\Omega$ and
since $\Omega$ is separating, we get $JW(\xi)J=W_{r}(J\xi)$.

\medskip

Recall that if $U_{t}$ is non trivial, the vacuum expectation
$\varphi$ is no longer tracial and is called the $q-$quasi-free
state. In fact in most cases (cf. \cite{H} Theorem 3.3), Araki-Woods
factors are type III von Neumann algebras.

\medskip

When $A$ is bounded, it is clear that our preliminaries are still valid with minor changes. For example we should get an extra $\|A^{-1}\|^{k/2}=\|A\|^{k/2}$ in the estimation of $\|U_{k}\|$. Note, in particular, that the Wick formula, as stated in Lemma 2, is still true, and that the following analogue of Bo\.zejko's scalar inequality holds : (proved in \cite{H})\\
 If $A$ is bounded, $(\eta_{u})_{u \in U}$ is a family of vectors in
$H^{\otimes n}$ and $(\alpha_{u})_{u \in U}$ a finitely supported
family of complex numbers then :
\begin{equation}
\left\|\sum_{u\in U}\alpha_{u} \eta_{u}\right\|_{q}
\le
\left\|\sum_{u\in U}\alpha_{u} W(\eta_{u})\right\|
\le
C_{|q|}^{\frac{3}{2}}
\frac{\|A\|^{\frac{n+1}{2}}-1}{\|A\|^{\frac{1}{2}}-1}\left\|\sum_{u\in U}\alpha_{u}
  \eta_{u}\right\|_{q}
\end{equation}
It is also a straightforward verification that Lemma 4, still hold in this setting. Observe also that $U$ is a unitary on $\FQH$ : this follows from the fact that for every $n \ge 1$, $P_{q}^{(n)}$, $A^{\otimes n}$ and $U$ commute on $H^{\otimes n}$. Note that $\mathcal{J}$ is no more an anti unitary from $H^{\otimes k}$ to $\overline{H^{\otimes k}}$, but since $U_{k}(I\otimes S)=M(a^{*}\otimes aU)$, we can deduce, as in the proof of Corollary 1, that $U_{k}(I\otimes S):\; H^{\otimes n-k}_{c}\otimes_{h} \overline{H^{\otimes k}_{r}}
\rightarrow B(\FT)$
is  completely bounded with norm less  than $C_{q}$, where $I$ stands for the identity of $H^{\otimes n-k}_{c}$. Following the same lines as in the proof of Theorem 1 we get :

\begin{thm}
Assume $A$ is bounded. Let $K$ be a complex Hilbert space. Then for all $n\ge 0$ and for all $\xi \in B(K)\otimes_{\min} H^{\otimes n}$ we have
\begin{equation}\label{uni}
\max_{0\le k \le n}\|( Id\otimes ((I\otimes S)R_{n,\;k}^{*})(\xi)\|
  \le \|({ Id} \otimes  W)(\xi)\|_{\min}
\end{equation}
$$ \phantom{\max_{0\le k \le n}\|( Id\otimes ((I\otimes S)R_{n,\;k}^{*})(\xi)\|
  \le}
\hspace{0.5 cm}\le C_{q}(n+1)
 \max_{0\le k \le n}\|({ Id}\otimes ((I\otimes S)R_{n,\;k}^{*})(\xi)\|$$
where ${ Id}$ denotes the identity mapping of $B(K)$, $I$ the identity of $H^{\otimes n-k}_{c}$, and where the
 norms of the left and right handsides are taken in
$B(K)\otimes_{\min}H^{\otimes n-k}_{c}\otimes_{\min}\overline{H^{\otimes k}_{r}}$.
\end{thm}

\bigskip

It is known (cf. \cite{H}) that if $U_{t}$ has a non trivial
continuous part then $\GQt$ is not injective. Using our techniques we are able to state a non-injectivity criterion similar to
that of  \cite{H} but independent of $q$.
\begin{cor} If either
$$\dim  E_{A}\left(\{1\}\right)H_{\C}\ge 2$$
or for some $T>1$
$$\frac{\dim E_{A}\left(]1,T]\right)H_{\C}}{T^{2}}> \frac{1}{2}$$
where $E_{A}$ is the spectral projection of $A$, then
  $\GQt$ is non injective.
\end{cor}

\proof
We can assume that $U_{t}$ is almost periodic, then we can write
$$(H_{\R},\; U_{t})=(\hat{H}_{\R},\;{\rm
  Id_{\hat{H}_{\R}}})\bigoplus_{\alpha \in
  \Lambda}(H_{\R}^{(\alpha)},\;U_{t}^{(\alpha)})$$
where
$$H_{\R}^{(\alpha)}=\R^{2},\;\;\;U_{t}^{(\alpha)}=
\left(
\begin{array}{cc}
\cos(t\ln \lambda_{\alpha}) & -\sin(t\ln \lambda_{\alpha})\\
\sin(t\ln \lambda_{\alpha}) & \cos(t\ln \lambda_{\alpha})
\end{array}
\right),
\;\;\; \lambda_{\alpha} >1$$
Thus the eigenvalues of the generator $A^{(\alpha)}$ of
$U_{t}^{(\alpha)}$ are $\lambda_{\alpha}$ and
$\lambda_{\alpha}^{-1}$.

\medskip

If $\dim E_{A}\left(\{1\}\right)H_{\C}\ge 2$ then $\dim \hat{H}_{\R} \ge 2$ and
since $U_{t}$ is trivial on $\hat{H}_{\R}$, the non-injectivity
follows from Theorem 2.

\medskip

For the remaining case we first suppose that $\dim H_{\R}=2$, $U_{t}$
is not trivial and  that $\GQt$ is injective. For all $n\ge 1$,
$A^{\otimes n}$ is a positive operator on $H^{\otimes n}$ equipped
with the deformed scalar product, we will denote by
$\lambda$ and $\lambda^{-1}$ the eigenvalues of $A$ with $\lambda >1$ and by $(\xi_{\un{i}})_{|\un{i}|=n}$ an orthonormal basis of eigenvectors of $A^{\otimes n}$ associated to the eigenvalues $(\lambda_{\un{i}})_{|\un{i}|=n}$. Since $\GQt$ is semidiscrete we must have for
every $n\ge 1$
\begin{equation*}
\|\sum_{|\un{i}|=n}W_{r}(J\xi_{\un{i}})W(\xi_{\un{i}})\|  \le
\|\sum_{|\un{i}|=n}W_{r}(J\xi_{\un{i}})\otimes W(\xi_{\un{i}})\| = \|\sum_{|\un{i}|=n}JW(\xi_{\un{i}})J\otimes
W(\xi_{\un{i}})\|
\end{equation*}
It is easily seen that
\begin{eqnarray*}
\|\sum_{|\un{i}|=n}W_{r}(J\xi_{\un{i}})W(\xi_{\un{i}})\| & \ge &
\sum_{|\un{i}|=n}\Dual{\Omega}{W_{r}(J\xi_{\un{i}})W(\xi_{\un{i}})\Omega}_{q}\\
& = & \sum_{|\un{i}|=n}\Dual{JW(\xi_{\un{i}})^{*}J\Omega}{W(\xi_{\un{i}})\Omega}_{q}\\
& = & \sum_{|\un{i}|=n}\Dual{\Delta^{\frac{1}{2}}\xi_{\un{i}}}{\xi_{\un{i}}}_{q} =
{\rm Trace}\left(\left(A^{-\frac{1}{2}}\right)^{\otimes n}\right)=(\lambda^{\frac{1}{2}}+\lambda^{-\frac{1}{2}})^{n}
\end{eqnarray*}
On the other hand, the map from $J\GQt J$ to $\overline{\GQt}$ taking $JW(\xi)J$ to $\overline{W(\xi)}$ is a $*$-isomorphism, hence
$$\|\sum_{|\un{i}|=n}JW(\xi_{\un{i}})J\otimes
W(\xi_{\un{i}})\|_{\min} = \|\sum_{|\un{i}|=n}\overline{W(\xi_{\un{i}})}\otimes
W(\xi_{\un{i}})\|_{\min}$$
Applying (\ref{uni}) twice, and recalling that on $H^{\otimes k}$,\; $S=J\Delta^{\frac{1}{2}}=J(A^{\otimes k})^{-\frac{1}{2}}$ and that $J :\;\overline{H_{r}^{\otimes k}}\rightarrow H_{r}^{\otimes k}$ is completely isometric, we get
$$
\|\sum_{|\un{i}|=n}\overline{W(\xi_{\un{i}})}\otimes
W(\xi_{\un{i}})\|_{\min} \le
C_{q}^{2}(n+1)^{2}
 \max_{0\le k,k^{\prime} \le n}\|\sum_{|\un{i}|=n}\overline{(I\otimes S)R_{n,\;k^{\prime}}^{*}(\xi_{\un{i}})}\otimes (I\otimes S)R_{n,\;k}^{*}(\xi_{\un{i}})\|$$
$$\hspace{1.9cm} \le C_{q}^{2}(n+1)^{2}
 \max_{0\le k,k^{\prime} \le n}
\|\sum_{|\un{i}|=n}\overline{(I\otimes (A^{\otimes k^{\prime}})^{-\frac{1}{2}})R_{n,\;k^{\prime}}^{*}(\xi_{\un{i}})}\otimes (I\otimes (A^{\otimes k})^{-\frac{1}{2}})R_{n,\;k}^{*}(\xi_{\un{i}})\|$$

Where the norms are computed in $\overline{H^{\otimes n-k^{\prime}}_{c}}\otimes_{\min}\overline{H^{\otimes k^{\prime}}_{r}}\otimes_{\min}H^{\otimes n-k}_{c}\otimes_{\min}H^{\otimes k}_{r}$. For a fixed $(k,k^{\prime})$, let us denote by
$$t = \sum_{|\un{i}|=n}\overline{(I\otimes (A^{\otimes k^{\prime}})^{-\frac{1}{2}})R_{n,\;k^{\prime}}^{*}(\xi_{\un{i}})}\otimes (I\otimes (A^{\otimes k})^{-\frac{1}{2}})R_{n,\;k}^{*}(\xi_{\un{i}})$$

As in the proof of Theorem 2, we have the following Hilbert-Schmidt estimate :
$$\|t\|_{\overline{H^{\otimes n-k^{\prime}}_{c}}\otimes_{\min}\overline{H^{\otimes k^{\prime}}_{r}}\otimes_{\min}H^{\otimes n-k}_{c}\otimes_{\min}H^{\otimes k}_{r}}\le
\|t\|_{\overline{H^{\otimes n-k^{\prime}}}\otimes_{2}\overline{H^{\otimes k^{\prime}}}\otimes_{2}H^{\otimes n-k}\otimes_{2}H^{\otimes k}}$$

\noindent Recall that $R_{n,\;k}^{*}:\;H^{\otimes n} \rightarrow H^{\otimes n-k}\otimes_{2}H^{\otimes k}$ is of norm less than $C_{|q|}^{\frac{1}{2}}$ and that $\|(A^{\otimes k})^{-\frac{1}{2}}\|_{B(H^{\otimes k})}=\lambda^{\frac{k}{2}}$. Hence,
\begin{eqnarray*}
\|t\|_{\overline{H^{\otimes n-k^{\prime}}}\otimes_{2}\overline{H^{\otimes k^{\prime}}}\otimes_{2}H^{\otimes n-k}\otimes_{2}H^{\otimes k}} & \le & C_{|q|}\lambda^{n}\|\sum_{|\un{i}|=n}\overline{\xi_{\un{i}}}\otimes \xi_{\un{i}}\|_{\overline{H^{\otimes n}}\otimes H^{\otimes n}}\\
& \le & C_{|q|}(\sqrt{2}\lambda)^{n}
\end{eqnarray*}

Combining all inequalities we get
$$(\lambda^{\frac{1}{2}}+\lambda^{-\frac{1}{2}})^{n}\le C_{|q|}^{3}(n+1)^{2}(\sqrt{2}\lambda)^{n}.$$

\medskip

We now return to the general case, we fix $T >1$ and we denote by
$\lambda_{1},\dots ,\lambda_{p}$ the eigenvalues of $A$ in $]1,T]$
counted with multiplicities. Thus we have $p=\dim
E_{A}\left(]1,T]\right)H_{\C}$. It is easy to deduce from our first step that for
any $n \ge 1$ we have
$$(\sum_{i=1}^{p}\lambda_{i}^{\frac{1}{2}}+\lambda_{i}^{-\frac{1}{2}})^{n}
\le C_{|q|}^{3}(n+1)^{2}(2p)^{\frac{n}{2}}T^{n}$$
Since for any $i$ we have
$\lambda_{i}^{\frac{1}{2}}+\lambda_{i}^{-\frac{1}{2}}\ge 2$ we deduce
$$(2p)^{n}\le C_{|q|}^{3}(n+1)^{2}(2p)^{\frac{n}{2}}T^{n}$$
So we necessarily have
$$\frac{2p}{T^{2}} \le 1$$
that is to say
$$\frac{\dim E_{A}]1,T]H_{\C}}{T^{2}}\le \frac{1}{2}$$
\Eproof

\bigskip
\noindent\underline{Acknowledgements} :\\
The author is grateful to Marek Bo\.zejko and Arthur Buchholz for
interesting discussions on these subjects and a very nice stay in
Wroclaw on July 2001, and to Eric Ricard for his careful reading and
improvements of the proofs.




\end{document}